\input amstex
\documentstyle{amsppt}
\magnification=1200
\NoBlackBoxes
\TagsOnRight
%
%
\define\scr{\Cal}

\define\s#1{{\frak #1}}

\define\e{{ \varepsilon }}
\define\oo{{ \infty }}
\define\floor#1{{ \lfloor {#1} \rfloor }}
\define\Ind#1{ 1_{ \{ #1 \} } }

\define\grad{\nabla}

\define\sym#1{{ #1^s }}
\define\compl#1{{ #1^c }}

\define\R{{\Bbb R}}

\define\Z{{\Bbb Z}}
\define\MM{{\Bbb M}}
\define\BB{{\Bbb B}}
\define\bdy{\partial}

\define\twoline#1#2{{\scriptstyle #1}\atop{\scriptstyle #2}}
\define\overbar{\overline}
\define\bull{\item"$\bullet$"}
\topmatter
\title From the P\'{o}lya-Szeg\"{o} symmetrization inequality for Dirichlet
integrals to comparison theorems for p.d.e.'s on manifolds$^*$\break
{\rm Presented at the 8th Symposium on Classical Analysis\break
Kazimierz Dolny, Poland, September, 1995}
\endtitle
\rightheadtext{symmetrization and p.d.e.'s on manifolds}
\author Alexander R. Pruss\endauthor
\affil Department of Philosophy\break
University of Pittsburgh\break
Pittsburgh, PA 15260\break
U.S.A.\break
e-mail: {\tt pruss+\@pitt.edu}
\endaffil
\thanks $^*$This research was partially supported by Professor
        J.\ J.\ F.\ Fournier's NSERC Grant \#4822 and was done
        and presented while
        the author was at the Department of Mathematics, University
        of British Columbia.
        The present notes are available on the Internet from
        {\it http://www.pitt.edu/$\sim$pruss/pss.html}\endthanks
\endtopmatter
\document
\remark{Disclaimer} These are only lecture notes.\endremark

\head 1. The Manifold Cases Under Consideration\endhead

Let $\MM_k^m$ be the simply connected constant curvature space form of
dimension $m$.
\roster
\bull $\MM_0^m$ is $\R^m$ with euclidean metric
\bull $\MM_k^m$ for $k>0$ is an $m$-sphere of radius $k^{-1/2}$
\bull $\MM_k^m$ for $k<0$ is $m$ dimensional hyperbolic space modelled on 
	the $m$-ball of radius $(-k)^{-1/2}$.
\endroster

Other definitions:
\roster
\bull all manifolds are Riemannian
\bull $M$ is an $m$-dimensional manifold, $m \ge 1$
\bull fix $k \in \R$
\bull for $\Omega \subseteq M$, let $\Omega^\#$ be the geodesic ball
	with the same volume as $\Omega$ centred about the origin $O$
	in $\MM_k^m$
\bull fix $B \subset M$, with infinitely differentiable boundary and
	homeomorphic to $\R^m$
\bull assume that for any open $\Omega \subseteq B$ with rectifiable boundary
we have the isoperimetric inequality
$$
	V_{m-1}(\bdy \Omega) \ge V_{m-1}(\bdy (\Omega^\#) ).
$$
\endroster

This isoperimetric inequality was conjectured by
Aubin~(1976) if all the sectional curvatures of $M$ are bounded above by $k$.
It is true if:
\roster
\bull $M=\MM_{k'}^m$ for $k' \le k$
\bull $m=2$ (Weil's first paper, 1926, etc.)
\bull $k\le 0$ and $m=3$ (Kleiner, 1992)
\bull $k=0$ and $m=4$ (Croke, 1984).
\endroster

Our results assume the inequality and only have substance where the
inequality is true.

For a real $f$ on a set $X$, let $f_t=\{ x : f(x)> t \}$.  Given a
map $(\cdot)^\#$ of measurable subsets of $X$ into measurable subsets of $Y$,
for $y \in Y$ let
$$
	f^\#(y) = \inf \{ t : y \in (f_t)^\# \}.
$$
If $(\cdot)^\#$ is measure preserving then $f^\#$ and $f$ are
equimeasurable.  If it is subset-preserving, then the Hardy-Littlewood 
inequality holds for real $f$ and $g$:
$$
	\int_X f g \le \int_Y f^\# g^\#.
$$

Let $f$ be non-negative and smooth on $M$, vanishing on $\bdy B$.
Then, it follows from the isoperimetric
inequality and coarea formula that
$$
	\int_B |\grad f|^2 \, dV_m \ge \int_{B^\#} |\grad f^\#|^2 \, dV_m.
$$
If $M=B=\MM_k^m=\R^m$, this is the P\'{o}lya-Szeg\"{o} Dirichlet-integral
symmetrization inequality.  Note that $\int_B |\grad f|^2=-\int_B f \Delta
f$.

Let $N$ be a manifold of dimension $n \ge 0$.  Given $\Omega \subseteq M
\times N$, let
$$
	\Omega^\# = \bigcup_{y \in N} (\Omega(y))^\# \times \{ y \},
$$
where $\Omega(y)=\{ x : (x,y) \in \Omega \}$.
The set $\Omega^\# \subseteq \MM_k^m \times N$ is a {\bf generalized Steiner
symmetrization}.
It is precisely Steiner symmetrization if $M=\R^m=\MM_k^m$ and $N=\R^n$.  The
operation $(\cdot)^\#$ induces a rearrangement on functions as above.

\head 2. Results on manifolds\endhead

\proclaim{Theorem}
Let $\Omega \subseteq B \times N$ have compact closure and nice boundary.
Let $u$ and $v$ be $C^2$ and non-negative on
$\Omega$ and $\Omega^\#$ respectively, vanishing outside their
respective domains, and
solving
$$
	-\Delta u = \phi(u)+\psi u+\lambda
$$
and
$$
	-\Delta v = \phi(v) + \psi^\# u + \lambda^\#,
$$
where $\phi$ is continuous decreasing on $[0,\oo)$ and $\lambda$ and $\psi$
continuous on $\bar B \times N$.  Let $\Phi$
be convex increasing on $[0,\oo)$.  Then, for all $y \in B$ we have
$$
	\int_{\Omega(y)} \Phi(u(x,y)) \, dV_m(x)
		\le \int_{\Omega(y)^\#} \Phi(v(x,y))\, dV_m(x),
$$
and moreover $v=v^\#$.
\endproclaim

In particular,
$$
	\max_x u(x,y) \le v(O,y),
$$
where $O$ is the origin in $\MM_k^m$.

A similar parabolic theorem can be proved, with the added condition that
$v$ satisfies the symmetrization of the initial condition for $u$.

\head 3. Some ideas for proofs\endhead

Given $u$ on $B \times N$, define on $B^\# \times N$:
$$
	u^I(x,y)=\int_{\BB(d(x,O))} u^\#(x',y) \, dV_m(x')
$$
(Baernstein $*$-function), where $\BB(r)$ is the ball around $O \in \MM_k^m$ with
geodesic radius $r$.
Given $v$ on $B^\# \times N$, define on $B^\# \times N$:
$$
	Jv(x,y)=\int_{\BB(d(x,O))} v(x',y) \, dV_m(x').
$$
Thus, $u^I=J(u^\#)$.
The conclusion of our theorem is equivalent to
$$
	u^I \le Jv.
$$

The proof of the theorem hinges on:
\proclaim{Proposition}	Let $u$ be as in the Theorem.  Then,
$$
	\int_{\Omega^\#} \vartheta \cdot (- \Delta u^\#)
		\le \int_{\Omega^\#} \vartheta \cdot ( \phi(u^\#) + \psi^\#
		u^\# + \lambda^\#),
$$
in the distributional sense
for every smooth (say $C^2$) function $\vartheta$ on $\Omega^\#$ vanishing on
the boundary and satisfying $\vartheta=\vartheta^\#$.
\endproclaim

Assume the Proposition.
Suppose we are in the case $\phi\equiv\psi\equiv 0$.  (The proof extends
to the general case by a clever method of Weitsman as in the case of
Steiner symmetrization on $\R^{n+m}$.)	Let $\vartheta$ be as in the
Proposition.  Then,
$$
	\int_{\Omega^\#} \vartheta \cdot (-\Delta v) = \int_{\Omega^\#}
	\vartheta \cdot \lambda^\#.
$$
By the Proposition if $g=u^\#-v$ then,
$$
	\int_{\Omega^\#} \vartheta \cdot (-\Delta g) \le 0.
$$
Let $\alpha$ be a continuous positive function on $\Omega^\#$, vanishing on
the boundary, and satisfying $\Delta \alpha = -1$ in $\Omega^\#$.  Let
$g_\e=g-\e \alpha$.  We have
$$
	\int_{\Omega^\#} \vartheta \cdot \Delta g_\e \ge \int_{\Omega^\#}
	\e \vartheta. $$
We shall show that it follows that $g_\e^I \le 0$ on $\Omega^\#$.  The
Theorem follows from this in the limit as $\e\to 0$.

Let $\scr M$ be the set of measures $\mu$ on
$\overbar{\Omega^\#}$ such that $\mu(A^\#) \ge \mu(A)$ for all Borel $A$ and
$\mu(\overbar{\Omega^\#}) \le 1$.  Define
$$
	G(\mu)=\int g_\e \, d\mu.
$$
We shall prove that $G \le 0$.	This will immediately imply that $g_\e^I(x,y)
\le 0$ since $g_\e^I(x,y)=|B(d(y,O))| G(\mu_{x,y})$ for an appropriate
probability measure $\mu_{x,y}$.  To prove $G \le 0$, let $\mu$ be the
measure at which $G$ attains a maximum, and assume that this maximum is
strictly positive.
Then $\mu$ has total mass $1$.
Taking slices carefully we may prove that there is an extremal
$\mu$ which has support contained in $M \times \{ y \}$ for some
$y$ and which is in be proportional to the measure $V_m$ on $M$
lifted to $M \times \{ y \}$ and restricted to some set $B(r) \times \{ y
\}$.  More precisely, for a continuous $f$
$$
	\int f d\mu = {1\over V_m(\BB(r))} \int_{V_m(\BB(r))} f(x,y) \, dV_m(x).
$$
I claim that the support of such $\mu$ cannot be contained inside
$\Omega^\#$.
For, if it is then for $t>0$ define a measure $\mu_t$ with density
$$
	\rho_t(z)=\int_{\Omega^\#} K_t(z,w) \, d\mu(w),
$$
on $\Omega^\#$, where $K_t$ is the heat kernel on $B^\#$ (Dirichlet boundary
conditions).  The measure $\mu_t$ will have total measure at most $1$.
Let $\mu_0=\mu$.
If the support of $\mu$ is contained in
$\Omega^\#$, then
$$
	\left. {d \over dt} \right|_{t=0+} \int_{\Omega^\#} g_\e d\mu_t
		= \int_{\Omega^\#} \Delta g_\e \, d\mu
		\ge \e.
$$
Of course this equation does not really make sense since $\Delta g_\e$ is only
defined distributionally, but we can make it make enough sense by
approximating $\mu$ with measures which have sufficiently smooth density.
The last inequality ``follows'' from the fact that $\int_{\Omega^\#}
\vartheta \Delta g_\e \ge \e$ if $\vartheta$ is positive, $C^2$, has
$\vartheta^\#=\vartheta$ and mean
$1$.  Since $\mu_t \in \scr M$ as it has desired symmetry because of the
symmetry of $K_t$, it follows that $\mu_0$ cannot be the extremal measure.

Now, suppose the extremal measure is proportional to a lifting of the measure
on $M$ to $\Omega^\#(y) \times \{ y \}$ for some $y$.
Let $r_0$ be such that $\Omega^\#(y)=\BB(r_0)$.  Assume that $r_0>0$.  (The
case $r_0=0$ is easy as $g_\e(0,y)=0$.)  Define $\nu_r$ to be
the measure on $B(r)$ lifted from the measure on $M$ so that
$$
	\int f d\nu_r = \int_{\BB(r)} f(x,y) \, d\nu_r(x).
$$
Then,
$$
	\left. {d \over dr} \right|_{r=r_0} \int g_\e \, d\nu_r
		= \left. { dV_m(\BB(r)) \over dr } \right|_{r=r_0}
			\cdot g_\e(x,y) = 0,
$$
where $x$ was such that $d(x,O)=r_0$.  Now then, let $\mu_r$ be $\nu_r$
normalized to have total mass $1$, i.e., let $\mu_r=(V_m(\BB(r)))^{-1} \nu_r$.
It follows that the derivative of
$$
	G(\mu_r) = {1\over V_m(\BB(r))} \int g_\e \, d\nu_r
$$
with respect to $r$ from the left at $r_0$ equals $G(\mu_{r_0})$ times the
derivative with respect to $r$ of $V_m(\BB(r))^{-1}$, which derivative is strictly
negative.  Thus, if $G(\mu_{r_0}) > 0$ then the derivative of $G(\mu_r)$ is
strictly negative at $r_0$, and it follows that $G(\mu_r)<G(\mu_{r'})$ for
some $r'<r$.  Since $\mu_{r'} \in \scr M$ for $r'\le r$, this is a
contradiction.	Hence $G(\mu_{r_0}) \le 0$, and the proof is complete.

To prove the Proposition, we first need the parabolic case of the Theorem with
$\phi=\lambda=0$.  This uses a slight extension of the Poly\'{a}-Szeg\"{o}
inequality, and
in effect has already been done by B\'{e}rard and Gallot (1980)\footnote{Note
added in 1997:  Cf.\ Gallot, 1988, Theorem~5.4(iii).}
This case can be rewritten as:
$$
	\int_{B^2} f(x) K_t^B(x,y) g(y) \, dV_{2m}(x,y)
		\le \int_{(B^\#)^2} f^\#(x) K_t^{B^\#}(x,y)
		g^\#(y) \, dV_{2m}(x,y),
$$
for $f,g \ge 0$ on $B$, where $K_t^B$ is the heat kernel vanishing on the
boundary of $B$.  This inequality is similar to the Riesz-Sobolev inequality.
Now, $K_t^{B \times N}((x_1,x_2),(y_1,y_2))=K_t^B(x_1,y_1) K_t^N(x_2,y_2)$.
This and Fubini's theorem implies that
$$
\eqalign{
	\int_{(B \times N)^2} &f(x) K_t^{B \times N}(x,y) g(y)
		\, dV_{2(m+n)}(x,y)\cr
		 &\le \int_{(B^\# \times N)^2} f^\#(x) K_t^{B^\# \times
		 N}(x,y) g^\#(y) \, dV_{2(m+n)}(x,y).}
$$

How can we use this to prove our Proposition?  Well, let $\vartheta$ and $u$
be as in
it.  Let $\tilde \vartheta$ be a function on $B \times N$ such that:
\roster
\bull $\tilde \vartheta^\#=\vartheta$
\bull $\tilde \vartheta$ is similarly ordered to $u$ (i.e.,
	$\tilde \vartheta(x) \le \tilde \vartheta(y)$ iff
	$u(x) \le u(y)$;  this is equivalent to requiring that
	$\int_{\Omega^\#} \tilde\vartheta^\#
	\cdot u^\# = \int_{\Omega} \tilde \vartheta \cdot u$.)
\endroster
We have
$$
	\eqalign{
		\int_{\Omega^\#} &\vartheta \cdot (- \Delta u^\#) \, dV_{m+n}
		\cr &=-\lim_{t\to 0+} {1\over t} \int_{(\Omega^\#)^2}
			[ \vartheta(x) K_t^{B^\# \times N}(x,y) u^\#(y)
			- \vartheta(x) \delta(x,y) u^\#(y) ] 
					\, dV_{2(m+n)}(x,y) \cr
		&\le -\lim_{t\to 0+} {1 \over t} \int_{\Omega^2}
			[ \tilde \vartheta(x) K_t^{B \times N}(x,y) u(y)
			- \tilde \vartheta(x) \delta(x,y) u(y) ] 
					\, dV_{2(m+n)}(x,y) \cr
		&= \int_\Omega \tilde \vartheta \cdot (- \Delta u) \,
		dV_{m+n} \cr
		&= \int_\Omega \tilde \vartheta \cdot ( \phi(u)+\psi u +
		\lambda)
		\cr &\le \int_{\Omega^\#} \vartheta \cdot (\phi(u^\#) +
		\psi^\# u^\#+\lambda^\#).
}
$$

\head 4. Discrete cases\endhead

The methods used can also give discrete symmetrization theorems.  Here, $M$
and $N$ are two discrete sets, and a laplacian is defined on $M \times N$.
Starting with a convolution-rearrangement inequality on $M$ like the one for
the heat kernel in the manifold case, one can duplicate most if not all of
the above work for a symmetrization on the discrete product set $M \times N$, 
with difference equations instead of p.d.e.'s.

An appropriate convolution-rearrangement inequality is known if $M$ is:
\roster
\bull the discrete line $\Z$, where we reorder functions so that
$$
	f^\#(0) \ge f^\#(1) \ge f^\#(-1) \ge \cdots
$$
(Hardy and Littlewood)
\bull the discrete circle $\Z_m$, where we reorder functions
so that
$$
	f^\#(0) \ge f^\#(1) \ge f^\#(-1) \ge \cdots
$$
	(effectively due to J.~R. Quine for the standard random
	walk;
	extended by the author to more general walks;  the limiting case as
	$m\to\oo$ is of course $\Z$)
\bull the $m$-regular tree $T_m$, with a spiral like reordering (this is
due to the author); this inequality also implies a Faber-Krahn
inequality for subsets of the $m$-regular tree
\bull the edge graph of an octahedron.
\endroster

All four convolution-rearrangement inequalities can be proved by a discrete
analogue of a method of Baernstein and Taylor~(1976), generalized (still in the
continuous case) by Beckner~(1993).

However, in the discrete case the method cannot handle many situations.
For instance, even the analogue of the Dirichlet integral inequality fails
on the graphs $\Z_2^3$ (cube) and $\Z_3^2$ (a euclidean plane based on a
finite field).\footnote{Note added in 1997: See Pruss, 1997a, 1997b, and
1997c for details of the discrete work.}

\head 5. Remarks added in 1997\endhead
It is worth noting that while the above symmetrization methods
symmetrize a manifold $M$ by using a manifold of revolution modelled on
the isoperimetric relations in $M$ (see Gallot, 1988) instead of
$\MM_k^n$.

The question of the general results that these kinds of methods can give
on manifolds is still open and the reader is invited to explore this
further.\footnote{The reader should feel free to contact the author if
the reader is interested in the project.} The present notes merely
outline the method. Further research on the manifold cases could
probably make use of the analogous but fully worked-out version of the
method in the discrete case (Pruss, 1997b; see especially Technical
Remark~3.1 for connections to manifolds).

\head Bibliography\endhead

Aubin (1976) [[ reference missing ]]

Baernstein II, A. and Taylor, B.A. (1976):
Spherical rearrangements, subharmonic functions and
        $*$-functions in $n$-space.  {\it Duke Math.\ J.}\ {\bf 43},
        245--268

Beckner, W. (1991):
Sharp {S}obolev inequalities on the sphere and the
        {M}oser-{T}rudinger inequality.
{\it Ann.\ Math.}\ {\bf 138}, 213--242

B\'{e}rard and Gallot (1980) [[ reference missing: instead see Gallot, 1988 ]]

Croke, C.B. (1984): A sharp four-dimensional
   isoperimetric inequality. {\it Comment.\ Math.\ Helv.}\ {\bf 59},
      187--192

Gallot, S. (1988):
In\'egalit\'es isop\'erim\'etriques et analytiques
     sur les vari\'et\'es riemanniennes.
     {\it Ast\'erisque} {\bf 163--164}, 31--91

Kleiner, B. (1992): An isoperimetric comparison theorem.
   {\it Invent.\ Math.}\ {\bf 108}, 37--47

Pruss, A.R. (1997a):
{Discrete convolution-symmetrization inequalities and
  the {F}aber-{K}rahn inequality on regular trees}.
  {\it Duke Math.\ J.}\ (to appear)

Pruss, A.R. (1997b):
       Symmetrization inequalities
       for difference equations on graphs. Preprint

Pruss, A.R. (1997c): {Discrete harmonic measure, {G}reen's functions and
  symmetrization: a unified probabilistic approach}. Preprint
\enddocument